\documentclass[12pt]{amsart}
\usepackage{amscd,amssymb,amsmath,latexsym}
\newtheorem{theorem}{Theorem}[section]
\newtheorem{lemma}[theorem]{Lemma}
\def\S{{\mathbb S}}
\def\Z{{\mathbb Z}}
\theoremstyle{definition}
\newtheorem{definition}[theorem]{Definition}
\newtheorem{example}[theorem]{Example}
\newtheorem{proposition}[theorem]{Proposition}
\newtheorem{corollary}[theorem]{Corollary}

\theoremstyle{remark}

\numberwithin{equation}{section}
\begin{document}
\title{On Fox spaces and Jacobi identities}
\author{Marek Golasi\'nski}
\address{Faculty of Mathematics and Computer Science, Nicolaus Copernicus University, Chopina 12/18, 87-100, Toru\'n, Poland}
\email{marek@mat.uni.torun.pl}
\author{Daciberg Gon\c calves}
\address{Dept. de Matem\'atica - IME - USP, Caixa Postal 66.281 - CEP 05311-970,
S\~ao Paulo - SP, Brasil}
\email{dlgoncal@ime.usp.br}
\author{Peter Wong}
\address{Department of Mathematics, Bates College, Lewiston, ME 04240, U.S.A.}
\email{pwong@bates.edu}
\thanks{This work was conducted during the second and third authors' visits to the Faculty of Mathematics and
Computer Science, Nicolaus Copernicus University, August 4 - 13, 2005 and July 29 - August 8, 2006. The second and third authors would
like to thank the Faculty of Mathematics and Computer Science for its hospitality and support. The third author was supported in part by a grant from Bates College and the National Science Foundation.}

\begin{abstract} In 1945, R.\ Fox introduced the so-called Fox torus homotopy
groups in which the usual homotopy groups are embedded and their
Whitehead products are expressed as commutators. A modern treatment of
Fox torus homotopy groups and their gene\-ralization has been given and studied. In this note, we further explore these groups and
their properties. We discuss co-multiplications on Fox spaces and
Jacobi identities for the generalized Whitehead products and the $\Gamma$-Whitehead products.
\end{abstract}
\date{May 3, 2007}
\maketitle

\newcommand{\af}{\alpha}
\newcommand{\et}{\eta}
\newcommand{\ga}{\gamma}
\newcommand{\ta}{\tau}
\newcommand{\ph}{\varphi}
\newcommand{\bt}{\beta}
\newcommand{\lb}{\lambda}
\newcommand{\wh}{\widehat}
\newcommand{\wt}{\widetilde}
\newcommand{\sg}{\sigma}
\newcommand{\om}{\omega}
\newcommand{\cH}{\mathcal H}
\newcommand{\cF}{\mathcal F}
\newcommand{\N}{\mathcal N}
\newcommand{\R}{\mathcal R}
\newcommand{\Ga}{\Gamma}
\newcommand{\F}{\mathfrak F}

\newcommand{\cc}{\mathcal C}

\newcommand{\bea} {\begin{eqnarray*}}
\newcommand{\beq} {\begin{equation}}
\newcommand{\bey} {\begin{eqnarray}}
\newcommand{\eea} {\end{eqnarray*}}
\newcommand{\eeq} {\end{equation}}
\newcommand{\eey} {\end{eqnarray}}

\newcommand{\ovl}{\overline}
\newcommand{\vv}{\vspace{4mm}}
\newcommand{\lra}{\longrightarrow}

\section*{Introduction}

In an attempt to give a geometric interpretation of the Whitehead product of \cite{whitehead}, R. Fox \cite{fox1,fox2} introduced the so-called torus homotopy groups $\tau_n(X,x_0)$ of a space $X$. For each $n\ge 1$, the group $\tau_n(X,x_0)$ is completely determined by the classical homotopy groups $\pi_k(X,x_0)$ for
$k\le n$ and the Whitehead products of the elements of $\pi_k(X,x_0)$. In fact, for any given $\alpha \in \pi_m(X,x_0)$ and $\beta\in \pi_n(X,x_0)$, the Whitehead product $[\alpha,\beta]$ when considered as an element in $\tau_k(X,x_0)$, $k\ge n+m-1$, is a commutator.

The Fox torus homotopy groups were given a modern interpretation and were generalized in \cite{ggw1}. Evaluation subgroups of $\tau_n$ were defined and their relationships with the classical Gottlieb groups were also studied. Moreover, it was shown in \cite{ggw1} that the generalized Whitehead product when embedded in the same spirit as \cite{arkowitz} in a larger
(and different) group is a commutator as well.

In 1954, M.\ Nakaoka and H.\ Toda \cite{nt} proved a Jacobi identity for
the classical Whitehead products involving elements in higher homotopy groups. In fact,
they gave two such proofs one of which made use of the Fox torus homotopy
groups. We should point out that H. Suzuki (1954)\cite{suzuki}, G. Whitehead (1954)\cite{gwhitehead2}, S.C. Chang (1954)\cite{chang}, P. Hilton (1955)\cite{hil}, H. Uehara and W. Massey (1957)\cite{uehara-massey}, and W. Barcus and M. Barratt (1958)\cite{barcus-barratt} obtained the same or similar Jacobi identities independently. Furthermore, Barcus and Barratt also derived a Jacobi identity in which one of the elements is in the fundamental group. In 1961, P.\ Hilton \cite{hilton} published a short note extending the Jacobi identity by allowing elements from the fundamental group.

\par For the generalized Whitehead product\footnote{The generalized Whitehead products
of \cite{ando,arkowitz,cohen} are slightly different from the others although they are all defined based upon commutators.}, the Jacobi identity was first proved by D. Cohen (1957)\cite[Section 5.8]{cohen}, in the language of carrier theory. It was pointed out by M. Arkowitz that such identity was also obtained in his doctoral dissertation (1960). When $A,B$ and $C$ are suspensions, the Jacobi identity was obtained independently by J. Boardman and B. Steer (1967)\cite[Equation (4.4)]{boardman-steer}, and by H. Ando (1968)\cite[Theorem 5.4]{ando}. The most general form of the Jacobi identity was due to
J.W.\ Rutter (1968)\cite{r,r1} who presented two proofs of the Jacobi identity for the generalized Whitehead product, under the assumption that $\Sigma A, \Sigma B$, and $\Sigma C$ are homotopy co-commutative. We should point out that J. Neisendorfer (1980) proved a Jacobi identity \cite[Proposition 9.10]{neisendorfer}, for the {\it external Samelson product} $\langle \cdot,\cdot\rangle$ for $[Z_1,G], [Z_2,G]$, and $[Z_3,G]$ where each $Z_i$ is {\it co-abelian} (i.e., $Z_i$ has weak category $\le 1$) and $G$ is group-like. In fact, $\langle \cdot,\cdot \rangle$ is the same as the generalized Whitehead product \cite[p.\ 467]{gwhitehead}. When $G=\Omega X$, Corollary 9.7\footnote
{When $Z$ is a connected $CW$-complex, this result was already obtained by I. Berstein and T. Ganea \cite[Corollary 3.3]{berstein-ganea}.} of \cite{neisendorfer} implies that $[Z,\Omega X]=[\Sigma Z,X]$ is an abelian group when $Z$ is co-abelian. This means that $\Sigma Z$ is homotopy co-commutative. Thus the Jacobi identity of Neisendorfer's follows from that of Rutter's. Furthermore, there are examples due to W. Gilbert \cite{gilbert} (see also \cite[Example 2.60, p. 70]{cornea})
where a space $X$ is not co-abelian while $\Sigma X$ is homotopy co-commutative.

The objective of this paper is to further explore the Fox torus homotopy groups and their generalizations. In particular, Section 1 takes up the study of Fox torus homotopy groups of spheres, Eilenberg-MacLane spaces,
and Fox spaces. Section 2 explores Fox spaces $F_n$ and their co-multiplication sets $\mathcal{C}(F_n)$
for $n\ge 1$. Section 3 is concerned with Jacobi identities for the generalized Whitehead product
studied in \cite{cohen,arkowitz,hilton3,boardman-steer,ando,neisendorfer,r,r1,oda}. By using the generalized Fox torus groups studied in \cite{ggw1},
we give a unified approach to the Jacobi identities of Rutter's \cite{r} and of Hilton's \cite{hilton}. In Section 4, we generalize the Fox torus groups to $\Gamma$-Fox torus groups for any
co-grouplike space $\Gamma$. We improve the Jacobi identity in \cite{oda} for the $\Gamma$-Whitehead product.
\par The authors are very grateful to Professor N.\ Oda for his very careful reading of the previous version of this
paper and all valuable comments and suggestions. We also thank Professor M. Arkowitz for bringing \cite{cornea} to our attention regarding spaces $X$ of weak category $>1$ while $\Sigma X$ is homotopy co-commutative.

\section{Fox spaces}

First, we recall from \cite{fox1,fox2} the definition of the $n$-th Fox torus homotopy group of
a pointed space $X$, for $n\ge 1$. Let $x_0$ be a basepoint of $X$, then
$$
\tau_n(X,x_0)=\pi_1(X^{T^{n-1}},\overline{x_0})
$$
where $X^{T^{n-1}}$ denotes the space of unbased maps from the $(n-1)$-torus $T^{n-1}$ to $X$ and
$\overline{x_0}$ is the constant map at $x_0$. When $n=1$, $\tau_1(X,x_0)=\pi_1(X,x_0)$.

To re-interpret Fox's result, we showed in \cite{ggw1} that
$$
\tau_n(X,x_0)\cong [\Sigma (T^{n-1}\sqcup *),X]
$$
the group of homotopy classes of basepoint preserving maps from the
reduced suspension of $T^{n-1}$ adjoined with a distinguished point
to $X$. As a result, we call $F_n =\Sigma (T^{n-1}\sqcup *)$ the $n$-th
{\em Fox space} with $F_1=\S^1$ the circle.

One of the main results of \cite{fox2} is the following split exact sequence:

\begin{equation}\label{fox-split}
0\to \prod_{i=2}^n \pi_i(X,x_0)^{\alpha_i(n)} \to \tau_n(X, x_0) \stackrel{\dashleftarrow}{\to} \tau_{n-1}(X, x_0) \to 1
\end{equation}
where $\alpha_i(n)$ is the binomial coefficient $\binom{n-2}{i-2}$.

With the isomorphism $\tau_{n-1}(\Omega X)\cong \prod_{i=2}^n \pi_i(X,x_0)^{\alpha_i(n)}$ shown in \cite[Theorem 1.1]{ggw1},
\eqref{fox-split} becomes

\begin{equation}\label{general-fox-split}
0\to \tau_{n-1}(\Omega X) \to \tau_n(X) \stackrel{\dashleftarrow}{\to} \tau_{n-1}(X) \to 1
\end{equation}
or equivalently
$$
0\to [\Sigma F_{n-1},X]\to [F_n,X]\stackrel{\dashleftarrow}{\to} [F_{n-1},X] \to 1.
$$

Using \eqref{fox-split}, the following is easy to verify.

\begin{proposition}\label{connectivity}
{\em For any pointed space $X$ with basepoint $x_0$, $\tau_k(X,x_0)=1$,
for all $k\le n$ if and only if $\pi_k(X,x_0)=1$, for all $k\le n$. Moreover,
$\tau_n(X,x_0)\cong \pi_n(X,x_0)$ if and only if $X$ is $(n-1)$-connected.}
\end{proposition}

For any $n\ge 1$, the $n$-th Fox space $F_n$
has the homotopy type of $\Sigma T^{n-1}\vee \S^1$. In fact, the natural isomorphism
$$[F_n,X]\cong[\Sigma F_{n-1},X]\rtimes[F_{n-1},X]$$ from the split exact sequence (\ref{general-fox-split}) yields:

\begin{proposition}\label{inductive-fox}
{\em For any integer $n> 1$, we have a homotopy equivalence
$$
F_n \simeq \Sigma F_{n-1} \vee F_{n-1}.$$}
\end{proposition}

As already pointed out in \cite{fox2}, one can define a Hurewicz homomorphism $\rho_n:\tau_n(X,x_0)\to H_n(X;\mathbb Z)$,
where $\Z$ is the group of integers. It follows from Proposition \ref{connectivity} that the (absolute) Hurewicz Isomorphism Theorem holds for the torus homotopy groups as in the classical case.

\begin{example}\label{ex1}
 Let $X=\S^n$, the $n$-sphere. For all $k \le n$, $\tau_k(\S^n)=\pi_k(\S^n)$. Using \eqref{fox-split}, we have $\tau_{n+1}(\S^n)\cong \pi_{n+1}(\S^n) \rtimes \pi_n(\S^n)$. More generally, we use \eqref{general-fox-split} to obtain
$$
\tau_m(\S^n)\cong \tau_{m-1}(\Omega \S^n) \rtimes \tau_{m-1}(\S^n).$$
\end{example}

\begin{example} \label{ex2}
Let $X=K(\pi,n)$ be an Eilenberg-MacLane space of type $(\pi,n)$ for some abelian group $\pi$ and $n\ge 2$. It follows from
\eqref{fox-split} that
\begin{equation}\label{semidirect}
\tau_m(X)\cong  \prod_{i=2}^m \pi_i(X,x_0)^{\alpha_i(m)} \rtimes \tau_{m-1}(X).
\end{equation}
Since $X=K(\pi,n)$, we get $\tau_m(X)\cong \pi_m(X)=1$ for $m<n$ and $\tau_n(X)=\pi_n(X)=\pi$. For $k\ge 1$, \eqref{semidirect} yields
$$
\tau_{n+k}(X) \cong \pi^{\alpha_n(n+k)} \rtimes \tau_{n+k-1}(X).
$$
The semi-direct product
structure comes from the Whitehead products which are all trivial since $X=K(\pi,n)$. Hence, we conclude that
$\tau_{n+k}(X)$ is a direct sum, i.e.,
$$
\tau_{n+k}(X) \cong \bigoplus_{j=0}^{k} {\pi_n(X)}^{\alpha_n(n+j)} \cong \pi^{\sum_{j=0}^k \binom{n-2+j}{n-2}}
$$
for any $k\ge 0$. Using a straighforward application of the Pascal's triangles one obtains that
$$
\sum_{j=0}^k\binom{n-2+j}{n-2}=\sum_{j=0}^k\binom{n-2+j}{j}=\binom{n+k-1}{k}.
$$
Therefore, for any space $X$ of type $K(\pi,n)$, we have
\[
   \tau_m(X)=
   \begin{cases}
      1,                       &\text{if $m< n$;}\\
      \pi^{\binom{m-1}{m-n}},  &\text{if $m\ge n$.}
   \end{cases}
\]
If $n=1$, that is, $X=K(\pi,1)$, then it is easy to show that $\tau_k(X)\cong \pi$, for all $k\ge 1$.
\end{example}
Next, we show that the wreath product $\wr$ appears in the second torus homotopy group of the Fox space $F_n$.
\begin{example}\label{ex3}
For any positive integer $n$, consider the $n$-th Fox space $F_n$. It was shown in \cite{ggw1} that $F_n$ has the homotopy type of a bouquet of spheres. More precisely,
$$F_n\simeq\bigvee_{i=1}^n {(\S^i)}^{\gamma_i(n)},$$
where $\gamma_i(n)$ is the binomial coefficient $\binom{n-1}{i-1}$. It follows that $\pi_1(F_n)\cong \mathbb Z$
and
$$\pi_2(F_n)\cong \pi_2(\tilde F_n)\cong \bigoplus_{j=-\infty}^{\infty} (\mathbb Z^{n-1})_j,$$
where $\tilde F_n$ denotes the universal cover of $F_n$. Now, \eqref{fox-split} leads
to the wreath product
$$
\tau_2(F_n)\cong \pi_2(F_n) \rtimes \tau_1(F_n) \cong \mathbb Z^{n-1} \wr \mathbb Z
$$
since the action of $\pi_1$ on the universal cover is the translation along the real line which covers the
only copy of $\S^1$ in $F_n$.
\end{example}

 Moreover, $F_m\simeq\bigvee_{i=1}^m {(\S^i)}^{\gamma_i(m)}$ is a homotopy retract of
 $F_n\simeq\bigvee_{i=1}^n {(\S^i)}^{\gamma_i(n)}$ provided that $m\le n$. Therefore, we derive:
\begin{proposition}\label{inclusion}
{\em Let $m\le n$. Then, there is an inclusion of groups
$$
\tau_k(F_m) \subseteq \tau_k(F_n)$$
for $k\ge 1$.}
\end{proposition}

\noindent{\bf Remark 1.1.}
The calculation of the torus homotopy groups of an Eilenberg-MacLane space $X$ of type $K(\pi,n)$ in Example
\ref{ex2} can be done directly as follows. Note that $F_m\simeq\bigvee_{i=1}^m {(\S^i)}^{\binom{m-1}{i-1}}$.
Therefore for $m\ge n$, the number of $n$-spheres in $F_m$ is precisely $\binom{m-1}{n-1}=\binom{m-1}{m-n}$.
The fact that the Whitehead products all vanish implies that $\tau_m(X)$ is isomorphic to the direct sum of $\binom{m-1}{m-n}$ copies of $\pi$. We should point out that the co-multiplication of $F_m$ does not yield a different group structure on $\tau_m(X)$ because $X$ has no non-trivial homotopy groups in dimensions other than $n$.

\section{Generalized Fox torus groups and co-multiplication}

In this section, we further explore the Fox spaces and examine
co-multiplications on them.
The functor $\tau_n$ was generalized in \cite{ggw1} as follows.

\begin{definition} \label{generalized-fox}
Let $X$ be a space and $x_0\in X$. For any space $W$, the $W$-{\em Fox group} of $X$ is defined to be
$$
\tau_W(X,x_0)=[\Sigma (W\sqcup *),X].
$$
\end{definition}

It is clear that $\tau_W$ reduces to $\tau_n$ when $W=T^{n-1}$.

We can further generalize Proposition \ref{inductive-fox} by the following construction.
For any pointed $W$, we let
$$
\F(W)= \Sigma W \vee W.
$$
Then, for any $n\ge 2$,
$$
\F_n(W) = \Sigma \F_{n-1}(W) \vee \F_{n-1}(W)
$$
or
$$
\F_n(W) = (\S^1 \wedge \F_{n-1}(W)) \vee \F_{n-1}(W)
$$
with $\F_1(W)=W$. Certainly, $\F_n(\S^1)\simeq F_n$. Furthermore, we can easily check
that
$$
[\F(W),X]=[W,A_{\F}(X)]
$$
for any pointed space $W$ where $A_{\F}(X)=X\times \Omega X$.

The definition of $\tau_W$ prompts another generalization of the Fox space as follows.

\begin{definition}
For any $n\ge 2$,
we let
$$
\hat {\F}_n(W)=\Sigma (W^{n-1} \sqcup *).
$$
\end{definition}

It follows that when $W=\S^1$, we have $F_n=\hat {\F}_n(\S^1) \simeq \F_n(\S^1)$.

Consider  the suspension co-multiplication on $\hat {\F}_n(\S^1)$.
While $\hat {\F}_n(\S^1)$ and $\F_n(\S^1)$ are homotopy equivalent, they are not
as co-$H$-spaces. To see that compare $[\hat {\F}_2(\S^1),F_2]$ with $[\F_2(\S^1),F_2]$.

Then, $[\hat {\F}_2(\S^1),F_2]=\tau_2(F_2) \cong
(\bigoplus_{j=-\infty}^{\infty}{(\mathbb Z)}_j) \rtimes \mathbb Z=\mathbb Z\wr\mathbb Z$, by the calculation in Example \ref{ex3}.
On the other hand, $\F_2(\S^1)=\Sigma \S^1 \vee \S^1=\S^2\vee \S^1$. It follows that
\begin{equation*}
\begin{aligned}
               {[\F_2(\S^1),F_2]} &\cong [\S^2 \vee \S^1, \S^2 \vee \S^1] \\
                &\cong [\S^2, \S^2 \vee \S^1] \oplus [\S^1, \S^2 \vee \S^1] \\
                &\cong (\bigoplus_{j=-\infty}^{\infty} {(\mathbb Z)}_j) \oplus \mathbb Z.
\end{aligned}
\end{equation*}
\par Following e.g., \cite[Proposition 1.1]{gg}, the kernel of the obvious epimorphism
$$[\F_n(\S^1),\F_n(\S^1)\vee \F_n(\S^1)]\to[\F_n(\S^1),\F_n(\S^1)\times\F_n(\S^1)]$$
is in 1-1 correspondence with the set $\mathcal{C}(\F_n(\S^1))$ of co-multiplications on the co-$H$-space $\F_n(\S^1)$
for $n\ge 1$. Because $\F_n(\S^1)\simeq F_n\simeq \bigvee_{i=1}^n(\S^i)^{\gamma_i(n)}$, we get the epimorphism
$$\big[\bigvee\nolimits_{i=1}^n(\S^i)^{\gamma_i(n)},\bigvee\nolimits_{i=1}^n(\S^i)^{\gamma_i(n)}\vee \bigvee\nolimits_{i=1}^n(\S^i)^{\gamma_i(n)}\big]\to
\big[\bigvee\nolimits_{i=1}^n(\S^i)^{\gamma_i(n)},\bigvee\nolimits_{i=1}^n(\S^i)^{\gamma_i(n)}\times\bigvee\nolimits_{i=1}^n(\S^i)^{\gamma_i(n)}\big].$$
Furthermore, $\pi_1(\F_n(\S^1))\cong\Z$, $\pi_k(\F_n(\S^1))\cong\pi_k(\bigvee_{j=-\infty}^\infty\bigvee_{i=2}^n(\S^i)^{\gamma_i(n)})$ and
$\pi_1(\F_n(\S^1)\vee\F_n(\S^1))\cong\Z*\Z$, $\pi_k(\F_n(\S^1)\vee\F_n(\S^1))\cong\pi_k(\bigvee_{j\in\Z*\Z}
\bigvee_{i=2}^n((\S^i)^{\gamma_i(n)}\vee(\S^i)^{\gamma_i(n)}))$ for $k\ge 2$, where $\mathbb Z\ast\mathbb Z$ is the free group on
two generators.
\par In particular, for $\F_1(\S^1)\simeq\S^1$ and $\F_2(\S^1)\simeq\S^2\vee \S^1$,
we get the epimorphisms
$$\Z*\Z\longrightarrow \Z\oplus \Z$$
and
$$\bigoplus_{j\in\Z\ast\Z}(\mathbb Z\oplus\mathbb Z)_j\oplus (\mathbb Z\ast\mathbb Z)\longrightarrow
\bigoplus_{j\in\Z\oplus\Z}(\mathbb Z\oplus\mathbb Z)_j\oplus (\mathbb Z\oplus\mathbb Z).$$

Hence, we can derive:
\begin{corollary} {\em There are bijections $\mathcal{C}(\F_1(\S^1))\cong (\mathbb Z\ast\mathbb Z)'$ and $\mathcal{C}(\F_2(\S^1))\cong (\mathbb Z\ast\mathbb Z)'$,
where $(\mathbb Z\ast\mathbb Z)'$ denotes the commutator subgroup of $\mathbb Z\ast\mathbb Z$ which is an infinitely
generated group.}
\end{corollary}

\noindent{\bf Remark 2.1.}
Observe that $\pi_k(\bigvee_{j\in J}\bigvee_{i=2}^n(\S^i)^{\gamma_i(n)})\cong\,
\mbox{colim}_{J'}\pi_k(\bigvee_{j'\in J'}\bigvee_{i=2}^n(\S^i)^{\gamma_i(n)})$
and $$\pi_k(\bigvee\nolimits_{j\in J}\bigvee\nolimits_{i=2}^n((\S^i)^{\gamma_i(n)}\vee(\S^i)^{\gamma_i(n)}))\cong\,\mbox{colim}_{J'}\pi_k(\bigvee\nolimits_{j'\in J'}
\bigvee\nolimits_{i=2}^n((\S^i)^{\gamma_i(n)}\vee(\S^i)^{\gamma_i(n)})),$$
where $J'\subseteq J$ runs over all finite subsets of $J$. Hence, we can make use of Hilton's result \cite{hil} so that the set of multiplications $\mathcal{C}(\F_n(\S^1))$ might be expressed by means of homotopy groups of appropriate spheres for $n\ge 1$.

\bigskip

Similar to the Fox space, the co-multiplication of $\hat{\F}_n(W)$ that yields
$\tau_{W^{n-1}}(X)=[\Sigma (W^{n-1} \sqcup *),X]=[\hat{\F}_n(W),X]$ is
given by the following split exact sequence which can be deduced from \cite[Theorem 3.1]{ggw1}:
\begin{equation}
1\to [\Sigma P_{n-1}(W),X] \to \tau_{W^{n-1}}(X) \stackrel{\dashleftarrow}{\to} \tau_{W^{n-2}}(X) \to 1
\end{equation}
which is the same as
$$
1\to [\Sigma P_{n-1}(W),X] \to [\hat{\F}_n(W),X]\stackrel{\dashleftarrow}{\to} [\hat{\F}_{n-1}(W),X] \to 1.
$$
Here, $P_n(W)=W^n/W^{n-1}$ is the $n$-fold pinched space of $W$.

Hence, $[\hat{\F}_n(W),X]\cong [\Sigma P_{n-1}(W),X]\rtimes[\hat{\F}_{n-1}(W),X]$ is
the semi-direct product with respect to a natural action
$$[\hat{\F}_{n-1}(W),X]\times[\Sigma P_{n-1}(W),X]\to[\Sigma P_{n-1}(W),X].$$
In particular, for $X=\hat{\F}_{n-1}(W)\vee\Sigma P_{n-1}(W)$, by the natural bijection
$[\hat{\F}_{n-1}(W),X]\times[\Sigma P_{n-1}(W),X]\cong[\hat{\F}_{n-1}(W)\vee\Sigma P_{n-1}(W),X]$,
the identity map $\mbox{id}_X$ is sent to the corresponding co-action
$$\alpha : \Sigma P_{n-1}(W)\to \hat{\F}_{n-1}(W)\vee \Sigma P_{n-1}(W).$$
Furthermore, the natural bijection $[\Sigma P_{n-1}(W),X]\rtimes[\hat{\F}_{n-1}(W),X]\cong
[\hat{\F}_n(W),X]$ for any space $X$ yields a homotopy equivalence
$$\hat{\F}_n(W)\simeq \Sigma P_{n-1}(W)\vee \hat{\F}_{n-1}(W).$$
Now, write $\nu : \Sigma P_{n-1}(W)\to \Sigma P_{n-1}(W)\vee \Sigma P_{n-1}(W)$
for the suspension co-multiplication on $\Sigma P_{n-1}(W)$ and
$\nu_{n-1} : \hat{\F}_{n-1}(W)\to\hat{\F}_{n-1}(W)\vee\hat{\F}_{n-1}(W)$
for the co-multiplication on $\hat{\F}_{n-1}(W)$.
\par Then, the maps
$$\nu^1_n : \hat{\F}_{n-1}(W)\stackrel{\nu_{n-1}}{\to}\hat{\F}_{n-1}(W)\vee\hat{\F}_{n-1}(W)\hookrightarrow \hat{\F}_n(W)\vee \hat{\F}_n(W)$$
and
{\small
\begin{equation*}
\nu^2_n : \Sigma P_{n-1}(W)\stackrel{\nu}{\to}\Sigma P_{n-1}(W)\vee \Sigma P_{n-1}(W)
\stackrel{\alpha\vee 1}{\to} \hat{\F}_{n-1}(W)\vee \Sigma P_{n-1}(W)\vee \Sigma P_{n-1}(W)
\hookrightarrow \hat{\F}_n(W)\vee \hat{\F}_n(W)
\end{equation*}
}

\noindent
lead to the suspension co-multiplication
$$\nu_n : \hat{\F}_n(W)\longrightarrow \hat{\F}_n(W)\vee \hat{\F}_n(W)$$ on the space
$\hat{\F}_n(W)$. Thus, we have shown:
\begin{theorem}\label{coH-structure}
For any $n> 1$, the co-multiplication on the homotopy type $\hat{\F}_n(W)\simeq \Sigma P_{n-1}(W)\vee\hat{\F}_{n-1}(W)$
of the space $\hat{\F}_n(W)$ that corresponds to the group structure on $\tau_{W^{n-1}}$ is given by the map constructed
above $$\nu_n : \hat{\F}_n(W)\longrightarrow \hat{\F}_n(W)\vee \hat{\F}_n(W).$$
\end{theorem}

Recall from \cite{ggw1} that the $n$-th Fox space $F_n=\hat{\F}^n(\S^1)$ has
the same homotopy type, as pointed spaces, of the pinched torus $T^n/T^{n-1}$
and by Proposition 2.2, $F_n\simeq\Sigma F_{n-1}\vee F_{n-1}$.
Hence, in light of Theorem \ref{coH-structure}, we derive:
\begin{corollary} {\em The co-multiplication on the homotopy type
$F_n\simeq\Sigma F_{n-1}\vee F_{n-1}$ of the Fox space $F_n$ that
corresponds to the group structure on $\ta_n$ is given by {\em Theorem \ref{coH-structure}} for $W=\mathbb S^1$.}
\end{corollary}
We end this section by noting that $P_n(W)$ is not homotopy equivalent to $\hat {\F^n}(W)$ in general.
For example, $P_2(T^2)=T^4/T^2$ whereas $\hat {\F^2}(T^2)=F_3$.

\section{Jacobi identities for generalized Whitehead products}

In \cite{hilton}, P. Hilton gave a Jacobi identity for the Whitehead product allowing elements of the fundamental group.
While the classical Jacobi identity (for elements in higher homotopy groups) of \cite{nt} has been generalized for the generalized Whitehead product, Hilton's result \cite{hilton} has not.
In this section, we give a unified approach to the Jacobi identities of Rutter's and Hilton's
using the generalized Fox torus groups studied in \cite {ggw1}. We extend Rutter's result by allowing $\Sigma C=\mathbb S^1$ or $\Sigma B=\Sigma C=\mathbb S^1$, thereby generalizing Hilton's result for the generalized Whitehead product. When one of the spaces is co-$H$, we can define the generalized Whitehead product via a group action similar to the treatment of Rutter's in \cite{r} except that Rutter assumed one of the spaces to be a suspension.

\par Now let $A$ and $B$ be pointed spaces.
Via the projections $p_A :((A\times B)\sqcup\ast) \to A$ and $p_B : ((A\times B)\sqcup\ast) \to B$,
we can regard $[\Sigma A,X]$ and $[\Sigma B,X]$ as subgroups of $\tau_{A\times B}(X)$.
Following \cite{hilton3}, recall that for any $\alpha \in [\Sigma A,X]$ and $\beta \in [\Sigma B,X]$, the generalized Whitehead product of $\alpha$ and $\beta$ is the element $[\alpha,\beta]=[K']$,
where $K':\Sigma (A\wedge B)\to X$ is induced by the map
$$
K=f'\cdot g'\cdot f'^{-1}\cdot g'^{-1}:\Sigma (A\times B) \to X
$$
which, when restricted to $\Sigma A \vee \Sigma B$, is homotopic to a constant.
Here, $f:\Sigma A \to X$, $g:\Sigma B \to X$ are maps representing $\alpha$ and
$\beta$, respectively and $f'=f\Sigma p_A, g'=g\Sigma p_B$.
Using the co-multiplication of $\Sigma (A\times B)$, $K$ is a well-defined map. It
is easy to observe that $[\alpha,\beta]=-(\Sigma t)^\ast[\beta,\alpha]$,
where $t :  A\wedge B\to B\wedge A$ is the switching equivalence of smash products.
\par In \cite{ggw1}, we gave another interpretation of the generalized
Whitehead product as follows\footnote{This differs from Theorem 4.1 of \cite{ggw1} in which the generalized Whitehead product of \cite{arkowitz} was used.}.
\begin{theorem}\label{whitehead}
Let $\alpha \in [\Sigma A, X]$ and $\beta \in [\Sigma B, X]$. Then,
the image of $[\alpha,\beta]$ in $\tau_{A\times B}(X)$ given by the
homotopy class of the composite
$$\Sigma((A\times B)\sqcup\ast)\to\Sigma(A\times B)\to\Sigma(A\wedge B)\stackrel{K'}{\to}X$$
is the commutator of the images of
$\alpha$ and $\beta$ in $\tau_{A\times B}(X)$.
\end{theorem}

For any pointed spaces $A$ and $B$ with disjoint unions $A \sqcup *_1$ and $B\sqcup *_2$, we have the following identification
$$
(A\sqcup *_1) \wedge (B \sqcup *_2)\approx (A\times B) \sqcup *_3,
$$
where $*_3$ corresponds to $(A\sqcup *_1) \vee (B \sqcup *_2)$ in the quotient $(A\sqcup *_1) \wedge (B \sqcup *_2)$. Since, by definition, $\tau_W(X)=[\Sigma (W\sqcup *),X]$, it follows that the generalized Whitehead product takes the form of
$$
\tau_A(X) \times \tau_B(X) \to \tau_{A\times B}(X).
$$
Hence, for $\alpha \in [\Sigma (A\sqcup *_1),X]$ and $\beta \in [\Sigma (B\sqcup *_2),X]$, the generalized Whitehead product $[\alpha, \beta]\in [\Sigma ((A\sqcup *_1) \wedge (B \sqcup *_2)),X]=\tau_{A\times B}(X)$ is a commutator.

Thus, the following diagram
\begin{equation}\label{fox-GWP}
\begin{CD}
    \tau_A(X) \times \tau_B(X) @>{(-,-)}>> \tau_{A\times B}(X)   \\
    @AAA                                 @AAA   \\
    [\Sigma A,X] \times [\Sigma B,X] @>{GWP}>> [\Sigma (A\wedge B),X]
\end{CD}
\end{equation}
is commutative.

By using the canonical maps $(A\sqcup *_1) \to A$ and $(B\sqcup *_2) \to B$, the image of $[\alpha,\beta]$ given by the composite map in Theorem \ref{whitehead} factors through $\tau_A(X) \times \tau_B(X)$. Furthermore, when $A=T^{m-1}$ and $B=T^{n-1}$, the generalized Whitehead product takes the form of
$$
GWP: \tau_m(X) \times \tau_n(X) \to \tau_{m+n-1}(X).
$$
This means that the generalized Whitehead product gives rise to a {\it Whitehead type product} on the Fox torus homotopy groups $\{\tau_k(X)\}$. Moreover,
if we denote by $WP$ the classical Whitehead product $\pi_m(X) \times \pi_n(X) \to \pi_{m+n-1}(X)$ then we have the following commutative diagram

\begin{equation}\label{fox-product}
\begin{CD}
    \tau_m(X) \times \tau_n(X) @>{GWP}>> \tau_{m+n-1}(X)   \\
    @A{j_m \times j_n}AA  @AA{j_{m+n-1}}A   \\
    \pi_m(X) \times \pi_n(X) @>{WP}>> \pi_{m+n-1}(X),
\end{CD}
\end{equation}

\noindent
where $j_k$ is the inclusion $\prod_{i=2}^k \pi_i(X,x_0)^{\alpha_i} \hookrightarrow \tau_k(X)$ in the split exact sequence \eqref{fox-split}
restricted to the only copy $\pi_k(X)$.

\par For any group $G$, Hall \cite[p.\ 150, 10.2.1.4]{hall} established the following {\em Jacobi identity} for
elements $a,b,c \in G$:
$$
(\star)\;\;(b^{-1},a,c)^b(c^{-1},b,a)^c(a^{-1},c,b)^a=1,
$$
where $(x,y,z)=((x,y),z)$, $(x,y)$ denotes the commutator $xyx^{-1}y^{-1}$
and $x^y=yxy^{-1}$ for $x,y,z\in G$. Furthermore, it can be easily shown
$$(\star\star)\;\;((y^{-1},x),z)^y=((x,y),z)((x,y),(z^{-1},y))^z.$$

To state our Jacobi identity for the generalized Whitehead product, we first recall from
\cite[Theorem 8.22]{bg} the following:

\begin{lemma} The natural transformations $A\wedge B\wedge C\to (A\wedge B)\wedge C$ and
$A\wedge B\wedge C\to A\wedge (B\wedge C)$ are homotopy equivalences provided that $A,B$ and $C$ have
the homotopy type of pointed compactly generated Hausdorff spaces.
\end{lemma}


\par Now, we prove the main result of this section. We shall write the operation in $\tau_{A\times B\times C}(X)$ additively.

\begin{theorem}\label{main}
Let $A,B$ and $C$ be pointed spaces with the homotopy type
of compactly generated Hausdorff spaces. Suppose $\alpha\in[\Sigma A,X]$, $\beta\in[\Sigma B,X]$ and $\gamma\in[\Sigma C,X]$ and denote by $\bar \alpha, \bar \beta,$ and $\bar \gamma$ the respective images in $\tau_{A\times B\times C}(X)$.
Then

\mbox{\em (1)}
$$(\Sigma t_{213})^*[[{\bar \beta}^{-1},{\bar \alpha}],{\bar \gamma}]^{\bar \beta}+(\Sigma t_{321})^*[[{\bar \gamma}^{-1},{\bar \beta}],{\bar \alpha}]^{\bar \gamma}+(\Sigma t_{132})^*[[{\bar \alpha}^{-1},{\bar \gamma}],{\bar \beta}]^{\bar \alpha}=1,$$
where $t_{ijk}$ is the appropriate twisting function for $i,j,k=1,2,3$.

\mbox{\em (2)}
If $\Sigma A$, $\Sigma B$ and $\Sigma C$ are homotopy co-commutative
co-$H$-spaces, then
$$[[\alpha,\beta],\gamma]+(\Sigma t_{312})^\ast[[\gamma,\alpha],\beta]+
(\Sigma t_{231})^\ast[[\beta,\gamma],\alpha]=0.$$

\mbox{\em (3)}
If each of $\Sigma A, \Sigma B,$ and $\Sigma C$ is either homotopy co-commutative or equal to $\mathbb S^1$ then
$$(\Sigma t_{213})^*{\beta}\cdot [[{\beta}^{-1},{\alpha}],{\gamma}]+(\Sigma t_{321})^*{\gamma}\cdot [[{\gamma}^{-1},{\beta}],{\alpha}]+(\Sigma t_{132})^*{\alpha}\cdot [[{\alpha}^{-1},{\gamma}],{\beta}]=1$$
where
\[
   \xi \cdot \eta=
                 \begin{cases}
                      \eta,            &\text{if $\xi \in [\Sigma W,X] \text{~and~} W\ne \mathbb S^0$;}\\
                      \eta^{\xi},      &\text{if $\xi \in \pi_1(X)$.}
                 \end{cases}
\]
\end{theorem}
{\bf Proof.} The first assertion (1) follows from the fact that the generalized Whitehead product $[-,-]$
in $\tau_{A\times B\times C}(X)$ coincides with the group theoretic commutator and the Jacobi identity
$(\star)$ together with the appropriate twisting functions.
\par For (2), we may assume, based upon the discussion above, that $\alpha,\beta,\gamma\in\tau_{A\times B\times C}(X)$.
Thence, by $(\star)$, we get
$$((\beta^{-1},\alpha),\gamma)^\beta
((\gamma^{-1},\beta),\alpha)^\gamma
((\alpha^{-1},\gamma),\beta)^\alpha=1$$
in the group $\tau_{A\times B\times C}(X)$.
But, by means of $(\star \star)$, $((\beta^{-1},\alpha),\gamma)^\beta=((\alpha,\beta),\gamma)((\alpha,\beta),(\gamma^{-1},\beta))^\gamma$.
Because of a homemorphism $\Sigma((A\times C)\wedge B)\approx(A\times C)\wedge\Sigma B$,
both elements $(\alpha,\beta)$ and $(\gamma^{-1},\beta)$ lie in the image of the abelian
group $[\Sigma((A\times C)\wedge B),X]$ in the group $\tau_{A\times B\times C}(X)$ and
we deduce that $((\alpha,\beta),(\gamma^{-1},\beta))^\gamma=1$.
Simi\-larly, $((\beta,\gamma),(\alpha^{-1},\gamma))^\alpha=
((\gamma,\alpha),(\beta^{-1},\alpha))^\beta=1$. Finally, in light
of Theorem \ref{whitehead},
the following equation $$[[\alpha,\beta],\gamma]+(\Sigma t_{231})^\ast[[\beta,\gamma],\alpha]+(\Sigma t_{312})^\ast
[[\gamma,\alpha],\beta]=0$$
holds in the image of the abelian group $[\Sigma(A\wedge B\wedge C),X]$ in the group
$\tau_{A\times B\times C}(X)$ and (2) is proven. This is the same as \cite[Theorem 1]{r} but our spaces need not be $CW$ complexes.

To prove (3), we note that the convention $\xi \cdot \eta$ is valid for the following reasons. When $\xi$ is not in $\pi_1(X)$, the image of $\eta^{\xi}$ is conjugation in the group $\tau$ and this will yield $\xi$ under the assumption that $\eta \in [\Sigma B, X]$ and $\Sigma B$ is co-commutative. On the other hand, if $\xi \in \pi_1(X)$ then $\eta^{\xi}$ is given by the classical action of $\pi_1(X)$. Now, (3) is simply the same as (1) when pulled back to $[\Sigma A\wedge B\wedge C,X]$. This is a generalization of formula (12) of Hilton's \cite{hilton}.\hfill$\square$

\bigskip

Note that the obvious pointed projection $A\sqcup\ast\to\mathbb{S}^0$ leads to the inclusion
$\pi_1(X)\subseteq \tau_A(X)$ for any nonempty topological space $A$ and a pointed space $X$.
In fact, it has been shown in \cite{ggw1} that $\tau_A(X)\cong [\Sigma A,X] \rtimes \pi_1(X)$
for any pointed spaces $A$ and $X$. Furthermore, this can be deduced from the following split exact
sequence of \cite[Theorem 3.1]{ggw1}
\begin{equation}\label{gen-split}
1\to [(V\times W)/V,\Omega X] \to \tau_{V\times W}(X) \stackrel{\dashleftarrow}{\to} \tau_V(X) \to 1
\end{equation}
which generalizes that of Fox.
\par When $V$ is a point and $W=A$, we have $\tau_A(X) \cong [\Sigma A,X] \rtimes \pi_1(X)$. The splitting \eqref{gen-split} gives the classical action of $\pi_1(X)$ on $[\Sigma A,X]$. Furthermore, when $V=A$ and $W=B$, this splitting gives rise
to an action of $[\Sigma A,X]$ ($\subseteq \tau_A(X)$) on $[(A\times B)/A,\Omega X]$. Since $[(A\times B)/A,\Omega X]=[\Sigma ((A\times B)/A), X]\cong [\Sigma (A\wedge B), X] \rtimes [\Sigma B,X]$, we have an action of
$[\Sigma A,X]$ on $[\Sigma (A\wedge B), X] \rtimes [\Sigma B,X]$.

Next, we give an alternative definition of the generalized Whitehead product when one of the spaces is co-$H$, thereby generalizing that of Rutter's \cite{r}.

When $B$ is a co-$H$ space, it follows from the proof of \cite[Theorem 3.1]{ggw1} that
$\Sigma ((A\times B)/A)\simeq \Sigma (B\vee (A\wedge B))$ is co-commutative so that $[\Sigma ((A\times B)/A), X]$ is abelian and $[\Sigma ((A\times B)/A), X]\cong [\Sigma (A\wedge B), X] \times [\Sigma B,X]$. Given $\alpha \in [\Sigma A,X]$ and $\beta \in [\Sigma B,X]$, we define
$$\alpha \ast \beta :=\alpha \cdot\tilde{\beta}$$ as the action of $\alpha$ on the image $\tilde{\beta}$  of $\beta$\
by the inclusion $[\Sigma B,X]\hookrightarrow[\Sigma (A\wedge B), X] \rtimes [\Sigma B,X]$.

Then, the image of the generalized Whitehead product $[\alpha,\beta]$ is given by the following relation in the group $\tau_{A\times B}(X)$:
\begin{equation}\label{action-GWP}
\overline{\alpha \ast \beta} =\overline{[\alpha,\beta]} + \bar \beta.
\end{equation}

Evidently, if $A$ is co-$H$, then one can define a similar action $\beta \ast \alpha$ so that \eqref{action-GWP} becomes $\overline{\beta \ast \alpha} =\overline{[\beta,\alpha]} + \bar \alpha$. These actions are analogous to the actions given in \cite[Theorem 2]{r}.

\section{$\Gamma$-Fox groups and Jacobi identities}

In \cite{oda}, N.\ Oda defined a generalization of the generalized Whitehead product ($\Gamma$-{\em Whitehead product})
by replacing the suspension $\Sigma$ with the smash product with a co-grouplike space $\Gamma$. A Jacobi identity was given in \cite[Theorem 2.11]{oda} for $[\Gamma \wedge A, X], [\Gamma \wedge B, X]$ and $[\Gamma \wedge C, X]$. In the case when $\Gamma=\mathbb S^1$, Oda established the same Jacobi identity as Rutter under the assumption that $A,B$ and $C$ be co-$H$ spaces. Thus, Rutter's result does not follow from Oda's since Oda's hypotheses imply that $\Sigma A, \Sigma B$ and $\Sigma C$ are necessarily homotopy co-commutative but a space $Z$ need not be co-$H$ while $\Sigma Z$ is homotopy co-commutative.

We end the paper by showing that the Fox group $\tau_W(X)$ can easily be generalized by replacing $\Sigma$ with $\Gamma$ and thus we can obtain the same kind of Jacobi identities as in the last section. These identites generalize those obtained in \cite{oda}.

\begin{definition} \label{gamma-fox}
Let $X$ be a space and $x_0\in X$ and $\Gamma$ a co-grouplike space. For any space $W$, the $W$-{\em $\Gamma$-Fox group} of $X$ is defined to be
$$
\tau^\Gamma_M(X)=\tau^{\Gamma}_W(X,x_0)=[\Gamma \wedge (W\sqcup *),X].
$$
\end{definition}
Note that
\begin{equation*}
      \tau^{\Gamma}_W(X)= [\Gamma (W\sqcup *),X]
                        = [(W\sqcup *) \wedge \Gamma,X]
                        = [(W\sqcup *), X^{\Gamma}_*],
\end{equation*}
where $X^{\Gamma}_*$ denotes the space of basepoint preserving maps from $\Gamma$ to $X$.
Following the notation of \cite{oda}, we write $\Gamma W=\Gamma \wedge W$.

\begin{theorem}\label{gamma-fox-split} Let $\Gamma$ be a co-grouplike space and $B$ a well-pointed space.
Then for any pointed space $X$ and any space $A$, the following sequence
\begin{equation}\label{split}
1\to [\Gamma ((A\times B)/A),X] \to \tau^{\Gamma}_{A\times B}(X) \stackrel{\dashleftarrow}{\to} \tau^{\Gamma}_A(X) \to 1
\end{equation}
is split exact.
\end{theorem}
\begin{proof}
The space $B$ is well-pointed, so we can consider the split cofibration
$$
A\sqcup * \stackrel{\dashleftarrow}{\to} (A\times B) \sqcup * \to (A\times B)/A.
$$
The corresponding Barrett-Puppe sequence yields the following short split-exact sequence
\begin{equation}\label{split-equation}
1\to [(A\times B)/A, Z] \to [(A\times B) \sqcup *, Z] \stackrel{\dashleftarrow}{\to} [A\sqcup *, Z] \to 1
\end{equation}
of pointed sets for any pointed space $Z$. For $Z=X^\Gamma_*$ using the adjoint isomorphism, \eqref{split-equation}
yields the desired short split-exact sequence.
\end{proof}

In particular, we have the following:

\begin{corollary}{\em  Let $W$ be a well-pointed space. Then the split cofibration $\ast\sqcup \ast\to W\sqcup\ast\to W$
leads to the isomorphism
$$\tau^{\Gamma}_W(X) \cong [\Gamma W,X] \rtimes [\Gamma, X].$$}
\end{corollary}

Given $\alpha \in [\Gamma A,X], \beta \in [\Gamma B,X]$, the $\Gamma$-{\it Whitehead product} $[\alpha,\beta]_{\Gamma}\in [\Gamma (A\wedge B),X]$ is defined (see \cite{oda}) in the same way as the generalized Whitehead product $GWP$ by simply replacing $\Sigma$ by $\Gamma$. The co-multiplication of $\Gamma$ allows the proof of Theorem 4.1 of \cite{ggw1} to remain valid when $GWP$ is replaced by $[-,-]_{\Gamma}$. Thus, the proof of Theorem 4.1 of \cite{ggw1} yields the following $\Gamma$-analog of Theorem \ref{whitehead}.

\begin{theorem}\label{gamma-whitehead}
Let $\alpha \in [\Gamma A, X]$ and $\beta \in [\Gamma B, X]$. Then,
the image of $[\alpha,\beta]_{\Gamma}$ in $\tau^{\Gamma}_{A\times B}(X)$ given by the
homotopy class of the composite
$$\Gamma((A\times B)\sqcup\ast)\to\Gamma(A\times B)\to\Gamma(A\wedge B)\stackrel{K'_{\Gamma}}{\to}X$$
is the commutator of the images of
$\alpha$ and $\beta$ in $\tau^{\Gamma}_{A\times B}(X)$.
\end{theorem}

Finally, we observe that the proof of Theorem \ref{main} remains valid when we replace $\Sigma$ with $\Gamma$
and the generalized Whitehead product with the $\Gamma$-Whitehead product. Thus, we have:

\begin{theorem}\label{gamma-main}
Let $A,B$ and $C$ be pointed spaces with the homotopy type
of compactly generated Hausdorff spaces and let $\Gamma$ be a co-grouplike space. Suppose
$\alpha\in[\Gamma A,X]$, $\beta\in[\Gamma B,X]$ and $\gamma\in[\Gamma C,X]$ and denote by $\bar \alpha, \bar \beta,$ and $\bar \gamma$ the respective images in $\tau^{\Gamma}_{A\times B\times C}(X)$.
Here $[-,-]_{\Gamma}$ denotes the corresponding $\Gamma$-Whitehead product.
Then

\mbox{\em (1)}
$$(\Gamma t_{213})^*[[{\bar \beta}^{-1},{\bar \alpha}]_{\Gamma},{\bar \gamma}]^{\bar \beta}_{\Gamma}+(\Gamma t_{321})^*[[{\bar \gamma}^{-1},{\bar \beta}]_{\Gamma},{\bar \alpha}]^{\bar \gamma}_{\Gamma}+(\Gamma t_{132})^*[[{\bar \alpha}^{-1},{\bar \gamma}]_{\Gamma},{\bar \beta}]^{\bar \alpha}_{\Gamma}=1.$$

\mbox{\em (2)}
If $\Gamma A$, $\Gamma B$ and $\Gamma C$ are homotopy co-commutative
co-$H$-spaces, then
$$[[\alpha,\beta]_{\Gamma},\gamma]_{\Gamma}+(\Gamma t_{312})^\ast[[\gamma,\alpha]_{\Gamma},\beta]_{\Gamma}+
(\Gamma t_{231})^\ast[[\beta,\gamma]_{\Gamma},\alpha]_{\Gamma}=0.$$

\mbox{\em (3)}
If each of $\Gamma A, \Gamma B$, and $\Gamma C$ is either homotopy co-commutative or equal to $\Gamma$ then
$$(\Gamma t_{213})^*{\beta}\cdot [[{\beta}^{-1},{\alpha}]_{\Gamma},{\gamma}]_{\Gamma}+(\Gamma t_{321})^*{\gamma}\cdot [[{\gamma}^{-1},{\beta}]_{\Gamma},{\alpha}]_{\Gamma}+(\Gamma t_{132})^*{\alpha}\cdot [[{\alpha}^{-1},{\gamma}]_{\Gamma},{\beta}]_{\Gamma}=1$$
where
\[
   \xi \cdot \eta=
                 \begin{cases}
                      \eta,            &\text{if $\xi \in [\Gamma W,X] \text{~and~} \Gamma W\ne \Gamma$;}\\
                      \eta^{\xi},      &\text{if $\xi \in [\Gamma, X]$.}
                 \end{cases}
\]
\end{theorem}

\noindent
{\bf Remark 4.1.} For the $\Gamma$-Whitehead product, one can give an alternative definition in terms of a group action similar to \eqref{action-GWP} under the hypothesis that one of the spaces is co-$H$.

\end{document}